\documentclass[12pt]{article}
\usepackage{hyperref}
\usepackage{amsfonts}
\usepackage{enumerate}
\usepackage{graphicx}
\usepackage{amsmath}
\begin{document}

\begin{center}
{\large\bf A boundary value problem of a generalised linear discrete time system with no solutions and infinitely many solutions}

\vskip.20in
Nicholas Apostolopoulos$^{1}$, Fernando Ortega$^{2}$ and\ Grigoris Kalogeropoulos$^{3}$\\[2mm]
{\footnotesize
$^{1}$National Technical University of Athens, Greece\\
$^{2}$ Universitat Autonoma de Barcelona, Spain\\
$^{3}$National and Kapodistrian University of Athens, Greece}
\end{center}

{\footnotesize
\noindent
\textbf{Abstract:}
In this article we study a class of generalised linear systems of difference equations with given boundary conditions and assume that the boundary value problem is non-consistent, i.e. it has infinite many or no solutions. We take into consideration the case that the coefficients are square constant matrices with the leading coefficient singular and provide optimal solutions. Numerical examples are given to justify our theory.
\\\\[3pt]
\\[3pt]

\vskip.2in

\section{Introduction}
Many authors have studied generalised discrete \& continuous time systems, see [1-28], and their applications, see [29-36]. Many of these results have already been extended to systems of differential \& difference equations with fractional operators, see [37-45]. We consider the generalised discrete time system of the form
\begin{equation}\label{eq1}
\begin{array}{cc}
FY_{k+1}=GY_k, & k=1,2,..., N-1
\end{array}
\end{equation}
with known boundary conditions
\begin{equation}\label{eq2}
\begin{array}{cc}A_1Y_0=B_1, &A_2Y_{N}=B_2.\end{array}
\end{equation}      
Where  $F, G \in \mathbb{R}^{r \times m}$, $Y_k\in \mathbb{R}^{m}$, and $A_1 \in \mathbb{C}^{r_1 \times m}$, $A_2 \in \mathbb{C}^{r_2 \times m}$, $B_1 \in \mathbb{C}^{r_1 \times 1}$ and $B_2 \in \mathbb{C}^{r_2 \times 1}$. The matrices $F$, $G$ can be non-square ($r\neq m$) or square ($r=m$) with $F$ singular (det$F$=0). 

Generalised linear systems of difference equations with given boundary conditions don't always guarantee to have unique solution. In the case where there exist solutions and they are infinite, we require optimal solutions for the system. The aim of this paper is to generalise existing results regarding the literature. An explicit and easily testable formula is derived of an optimal solution for the system. 

%The paper is organized as follows: in Section 2 we 
%refer to the mathematical background used throughout this paper, in Section 3 we study the optimal solutions of system (1) with given non-consistent initial conditions of type (2) in Section 4 we study the cases of a non-consistent boundary value problem in the form of (3)-(4) and of a consistent with infinite solutions and in Section 5 we provide a reformulated Kalman filter for singular non-homogeneous linear control systems of fractional nabla difference equations.

\section{Preliminaries}
Throughout the paper we will use in several parts matrix pencil theory to establish our results. A matrix pencil is a family of matrices $sF-G$, parametrized by a complex number $s$, see [46-53].
\\\\
\textbf{Definition 2.1.} Given $F,G\in \mathbb{R}^{r \times m}$ and an arbitrary $s\in\mathbb{C}$, the matrix pencil $sF-G$ is called:
\begin{enumerate}
\item Regular when  $r=m$ and  det$(sF-G)\neq 0$;
\item Singular when  $r\neq m$ or  $r=m$ and det$(sF-G)\equiv 0$.
\end{enumerate}
In this article we consider the system \eqref{eq1} with a \textsl{regular pencil}, where the class of $sF-G$ is characterized by a uniquely defined element, known as the Weierstrass canonical form, see [46-53], specified by the complete set of invariants of $sF-G$. This is the set of elementary divisors of type  $(s-a_j)^{p_j}$, called \emph{finite elementary divisors}, where $a_j$ is a finite eigenvalue of algebraic multiplicity $p_j$ ($1\leq j \leq \nu$), and the set of elementary divisors of type $\hat{s}^q=\frac{1}{s^q}$, called \emph{infinite elementary divisors}, where $q$ is the algebraic multiplicity of the infinite eigenvalue. $\sum_{j =1}^\nu p_j  = p$ and $p+q=m$.
\\\\
From the regularity of $sF-G$, there exist non-singular matrices $P$, $Q$ $\in \mathbb{R}^{m \times m}$ such that 
\begin{equation}\label{eq3}
\begin{array}{c}PFQ=\left[\begin{array}{cc} I_p&0_{p,q}\\0_{q,p}&H_q\end{array}\right],
\\\\
PGQ=\left[\begin{array}{cc} J_p&0_{p,q}\\0_{q,p}&I_q\end{array}\right].\end{array}
\end{equation}
$J_p$, $H_q$ are appropriate matrices with $H_q$ a nilpotent matrix with index $q_*$, $J_p$ a Jordan matrix and $p+q=m$. With $0_{q,p}$ we denote the zero matrix of $q\times p$. The matrix $Q$ can be written as
\begin{equation}\label{eq4}
Q=\left[\begin{array}{cc}Q_p & Q_q\end{array}\right].
\end{equation}
$Q_p\in \mathbb{R}^{m \times p}$ and $Q_q\in \mathbb{R}^{m \times q}$. The matrix $P$ can be written as
\begin{equation}\label{eq5}
P=\left[\begin{array}{c}P_1 \\ P_2\end{array}\right].
\end{equation}
$P_1\in \mathbb{R}^{p \times r}$ and $P_2\in \mathbb{R}^{q \times r}$.
%\\\\
%\textbf{Definition 2.2.}  Let $J_p$ be a Jordan matrix as defined in (2). Then with $F_{n,n}(J_p(k+n)^{\bar n})$ we will denote the discrete Mittag-Leffler function with two parameters defined by 
%\begin{equation}
%F_{n,n}(J_p(k+n)^{\bar n})=\sum^{\infty}_{i=0}J_p^i\frac{(k+n)^{\overline{in}}}{\Gamma((i+1)n)}.
%\end{equation}
The following results have been proved.
\\\\
\textbf{Theorem 2.1.}  (See [1-28]) We consider the systems \eqref{eq1} with a regular pencil. Then, its solution exists and for $k\geq 0$, is given by the formula
\begin{equation}\label{eq6}
    Y_k=Q_pJ_p^kC.  
\end{equation}
Where $C\in\mathbb{R}^p$ is a constant vector. The matrices $Q_p$, $Q_q$, $P_1$, $P_2$, $J_p$, $H_q$ are defined by \eqref{eq3}, \eqref{eq4}, \eqref{eq5}. 
\\\\
A generalised boundary value problem (BVP) of type \eqref{eq1}--\eqref{eq2} can have a unique solution, infinite solutions or no solutions, see [8], [15].  
\\\\
\textbf{Definition 2.2.} If for the system \eqref{eq1} with boundary conditions \eqref{eq2} there exists at least one solution, the BVP \eqref{eq1}--\eqref{eq2} is said to be consistent.
\\\\
Let
\begin{equation}\label{eqK}
K=\left[\begin{array}{c} A_1Q_p\\A_2Q_p(N+1)^{\overline{n-1}}F_{n,n}(J_p(N+n)^{\bar n})(I_p-J_p) \end{array}\right]
\end{equation}
and
\begin{equation}\label{eqL}
L=\left[\begin{array}{c} B_1\\B_2 \end{array}\right].
\end{equation}
Then
\\\\
\textbf{Theorem 2.2.} (See [8], [15]) We consider the BVP \eqref{eq1}--\eqref{eq2} with a regular pencil.
Then the BVP is consistent,  if and only if one of the following is satisfied:
\begin{itemize} 
\item $p= r_1+r_2$ and rank$K=p$;
\item $p< r_1+r_2$, rank$K=p$ and $L\in colspan K$;
\item $p> r_1+r_2$, rank$K=r_1+r_2$;
\item $L\in colspanK$ and $K$ is rank deficient.
\end{itemize}

\section{Non-consistent boundary value problem}

We can now state the following Theorem, which provides optimal solutions for non-consistent BVPs of type \eqref{eq1}--\eqref{eq2}.
\\\\
\textbf{Theorem 3.1.} We consider the system \eqref{eq1} with known boundary conditions of type \eqref{eq2} with a regular pencil.
Then for the  non-consistent BVP of type \eqref{eq1}--\eqref{eq2}
\begin{enumerate}[(a)]
\item If $p< r_1+r_2$, $L\notin colspanK$ and $K$ is full rank, an optimal solution of the BVP is given by
\begin{equation}\label{eqT1}
    \hat Y_k=Q_pJ_p^k(K^*K)^{-1}K^*L.
\end{equation}
\item If $L\notin colspanK$ and $K$ is rank deficient, an optimal solution of the BVP is given by
\begin{equation}\label{eqT2}
    \hat Y_k=Q_pJ_p^k(K^*K+E^*E)^{-1}K^*L.
\end{equation}
Where $E$ is a matrix such that $K^*K+E^*E$ is invertible and $\left\|E\right\|_2=\theta$, $0<\theta<<1$.
\end{enumerate}

The matrices $K$, $L$ are given by \eqref{eqK}--\eqref{eqL} and the matrices $Q_p$, $J_p$ are given by \eqref{eq3},\eqref{eq4}.\\\\
\textbf{Proof.} From Theorem 2.1 it is known that the solution of the system \eqref{eq1} exists and is given by \eqref{eq6}. For given boundary conditions of type \eqref{eq2}, $C$ is the solution of the linear systems
\[
A_1Q_pC=B_1
\]
and
\[
A_2Q_p(N+1)^{\overline{n-1}}F_{n,n}(J_p(N+n)^{\bar n})(I_p-J_p)C=B_2,
\]
or, equivalently,
\begin{equation}\label{eqKC}
KC=L.
\end{equation}
System \eqref{eqKC} has $r_1+r_2$ linear equations and $p$ unknowns. 
\\\\
For the the proof of (a), since $p< r_1+r_2$, $L\notin colspanK$ and $K$ is full rank, the system \eqref{eqKC} has no solutions. Let $\hat L$ be a vector such that a vector $\hat C$ is the unique solution of the system $K\hat C=\hat L$. We want then to solve the following optimisation problem
\[
\begin{array}{c}min\left\|L-\hat L\right\|_2^2\\s.t.\quad K\hat C=\hat L,\end{array}
\]
or, equivalently,
\[
min\left\|L-K\hat C\right\|_2^2.
\]
Where $\hat C$ is the optimal solution, in terms of least squares, of the linear system \eqref{eqKC}. In this case, with similar steps as [38], [54], [55], the solution $\hat C$ is given by
\[
\hat C=(K^*K)^{-1}K^*L.
\]
This is the least squares solution and hence an optimal solution of the system \eqref{eq1} with boundary conditions of type \eqref{eq2} is given by \eqref{eqT1}. 
\\\\
For the proof of (b), since $L\notin colspanK$ and $K$ is rank deficient, the system \eqref{eqKC} has no solutions. While the matrix $K$ is rank deficient, the matrix $K^*K$ is singular and hence not invertible. Thus we can not apply the same method as in (a). In this case, we seek a solution $\hat C$ minimizing the functional
\[
D_2(\hat C)=\left\|L-K\hat C\right\|_2^2+\left\|E\hat C\right\|_2^2.
\]
Where $E$ is a matrix such that $K^*K+E^*E$ is invertible and $\left\|E\right\|_2=\theta$, $0<\theta<<1$. Expanding $D_2(\hat C)$ gives
\[
D_2(\hat C)=(L-K\hat C)^*(L-K\hat C)+(E\hat C)^*E\hat C,
\]
or, equivalently,
\[
D_2(\hat C)=L^*L-2L^*K\hat C+(\hat C)^*K^*K\hat C+(\hat C)^*E^*E\hat C,
\]
because $L^*K\hat C=(\hat C)^*K^*L$. Furthermore
\[
\frac{\partial}{\partial \hat C }D_2(\hat C)=-2K^*L+2K^*K\hat C+2E^*E\hat C.
\]
Setting the derivative to zero, $\frac{\partial}{\partial \hat C }D_2(\hat C)=0$, we get
\[
(K^*K+E^*E)\hat C=K^*L.
\]
The solution is then given by
\[
\hat C=(K^*K+E^*E)^{-1}K^*L.
\]
Hence an optimal solution in this case is given by \eqref{eqT2}. The proof is completed.
\\\\
\textbf{Remark 3.1.} In Theorem 3.1 the optimal solution of the BVP problem \eqref{eq1}--\eqref{eq2} was achieved after a perturbation to the column $L$ accordingly 
\[
min\left\|L-\hat L\right\|_2,
\]
or, equivalently,
\[
\left\|L-K(K^*K)^{-1}K^*L\right\|_2,
\]
in case of (a) and 
\[
\left\|L-K(K^*K+E^*E)^{-1}K^*L\right\|_2,
\]
in case of (b). 
\\\\
A consistent BVP doesn't guaranty unique solution. The case of existence and uniqueness of solutions has been fully analysed in [8], [15]. 
\\\\
\textbf{Proposition 3.1.} (See [8], [15]) When the BVP \eqref{eq1}--\eqref{eq2} is consistent, it has a unique solution if and only if one of the following is satisfied
\begin{itemize} 
\item $p= r_1+r_2$ and rank$K=p$;
\item $p< r_1+r_2$, rank$K=p$ and $L\in colspan K$.
\end{itemize}
For a consistent BVPs of type \eqref{eq1}--\eqref{eq2} with infinite solutions, there are different ways to obtain optimal solutions, depending on the problem that the BVP represents. In the following Proposition, we apply some techniques in order to gain the minimum solution.
\\\\
\textbf{Proposition 3.2.} We consider the system \eqref{eq1} with known boundary conditions of type \eqref{eq2} and a regular pencil. Then in both cases, for a consistent BVP of type \eqref{eq1}--\eqref{eq2} with infinite solutions:
\begin{enumerate}[(a)]
\item If $K^\dagger$ is the Moore-Penrose Pseudoinverse of $K$, then an optimal solution of the BVP \eqref{eq1}--\eqref{eq2} is given by 
\begin{equation}\label{eqP1}
    \hat Y_k=Q_pJ_p^kK^\dagger L.
\end{equation}\label{eqP2}
\item If $p> r_1+r_2$, $K$ is full rank, a minimum solution of the BVP is given by
\begin{equation}\label{eqP3}
    \hat Y_k=Q_pJ_p^k(KK^*)^{-1}KL.
\end{equation}
\item If $L\in colspanK$ and $K$ is rank deficient, an optimal solution of the BVP is given by \eqref{eqT2}.
\end{enumerate}
The matrices $K$, $L$ are given by \eqref{eqK}, \eqref{eqL} and the matrices $Q_p$, $J_p$ are given by \eqref{eq3},\eqref{eq4}.\\\\
\textbf{Proof.} From Theorem 2.1 it is known that the solution of the system \eqref{eq1} exists and is given by \eqref{eq6}. For given boundary conditions of type \eqref{eq2}, $C$ is the solution of the linear system \eqref{eqKC} which has $r_1+r_2$ linear equations and $p$ unknowns. 
\\\\
For the proof of (a), given an $(r_1+r_2)\times p$ matrix $K$, the Moore-Penrose Pseudoinverse of $K^\dagger$ is calculated from the singular value decomposition of $K$, see [19, 26]. The solution of system \eqref{eqKC} is then given by $\hat C=K^\dagger L$.
\\\\
For the proof of (b), $p> r_1+r_2$, $K$ is a wide matrix (more columns than rows) and full rank, i.e. the system \eqref{eqKC} is underdetermined. In this case, it is common to seek a solution $\hat C$ (one of the infinite solutions of system \eqref{eqKC}) with minimum norm. Hence we would like to solve the optimisation problem
\[
\begin{array}{c}\quad min \quad \left\|\hat C\right\|_2^2,\\
\quad s.t.\quad L=K\hat C.
\end{array}
\]
By defining the Lagrangian
\[
D_3(\hat C,\lambda)=\left\|\hat C\right\|_2^2+\lambda^*(L-K\hat C).
\]
and taking the derivatives of the Lagrangian we get
\[
\frac{\partial}{\partial \hat C }D_3(\hat C,\lambda)=2\hat C- K^*\lambda
\]
and
\[
\frac{\partial}{\partial \lambda }D_3(\hat C,\lambda)=L-K\hat C.
\]
Setting the derivatives to zero, $\frac{\partial}{\partial \hat C }D_3(\hat C,\lambda)=0$ and $\frac{\partial}{\partial \lambda }D_3(\hat C,\lambda)=0$, we get
\[
\hat C=\frac{1}{2} K^*\lambda
\]
and
\[
L=K\hat C,
\]
or, equivalently,
\[
L=\frac{1}{2}KK^*\lambda.
\]
Since $rank K= r_1+r_2$, the matrix $KK^*$ is invertible and thus
\[
\lambda=2(KK^*)^{-1}L.
\]
The solution is then given by
\[
\hat C=K^*(KK^*)^{-1}L.
\]
Hence an optimal solution in this case is given by \eqref{eqP2}.
\\\\
For the proof of (c), since $L\in colspanK$ and $K$ is rank deficient, the system \eqref{eqKC} has infinite solutions. While matrix $K$ is rank deficient, the matrix $KK^*$ is singular and hence not invertible. In this case, we seek a solution $\hat C$ minimizing the functional
\[
D_4(\hat C)=\left\|L-K\hat C\right\|_2^2+\left\|E\hat C\right\|_2^2
\]
and by working as in Theorem 3.1 (b), the solution is given by
\[
\hat C=(K^*K+E^*E)^{-1}K^*L
\]
and hence an optimal solution in this case is given by \eqref{eqT2}. The proof is completed.

\subsection*{Numerical examples}

\textit{Example 4.1.}
\\\\
Assume the system \eqref{eq1} for $k$=0, 1, 2, 3. Let
\[
F=\left[\begin{array}{ccccc} 0&1&0&0&0\\0&0&0&0&-1\\1&-1&0&1&-1\\1&-1&0&2&-1\\0&0&0&0&0\end{array}\right]
\]
and
\[
G=\left[\begin{array}{ccccc} 0&0&1&-1&1\\0&0&0&0&-\frac{1}{4}\\\frac{1}{2}&-\frac{1}{2}&0&-\frac{1}{2}&\frac{1}{2}\\\frac{1}{2}&-\frac{1}{2}&1&-\frac{1}{2}&\frac{3}{4}\\0&0&0&1&-1\end{array}\right].
\]
Then det($sF-G$)=$(s-\frac{1}{2})s(s-\frac{1}{4})$ and the pencil is regular. We assume the boundary conditions \eqref{eq2} with 
\[
A_1=\left[\begin{array}{ccccc}1&0&0&0&0\\0&0&0&1&-1\\0&0&1&0&0\\0&0&1&0&0\end{array}\right],
\]
\[
A_2=\left[\begin{array}{ccccc}0&0&0&0&1\end{array}\right],
\]
\[
B_1=\left[\begin{array}{c}0\\0\\36\\0\end{array}\right],
\]
and
\[
B_2=24
\]
or, equivalently,
\[
L=\left[\begin{array}{c}0\\0\\36\\0\\24\end{array}\right].
\]
The three finite eigenvalues ($p$=3) of the pencil are $\frac{1}{2}$, 0, $\frac{1}{4}$ and the Jordan matrix $J_p$ has the form
\[
J_p=\left[\begin{array}{ccc}\frac{1}{2}&0&0\\0&0&0\\0&0&\frac{1}{4}\end{array}\right].
\]
By calculating the eigenvectors of the finite eigenvalues we get the matrix $Q_p$
\[
Q_p=\left[\begin{array}{ccc}1&1&0\\0&1&0\\0&0&0\\0&0&1\\0&0&1\end{array}\right].
\]
After computations
\[
K=\left[\begin{array}{ccc}1&1&0\\0&0&0\\0&0&0\\0&0&0\\0&0&\frac{24}{36}\end{array}\right]
\]
Furthermore $r_1+r_2=5>3=p$, $rank K=2<3=p$ and $L \notin colspan K$, i.e. the BVP is non-consistent. Thus from Theorem 3.1 and (b), an optimal solution of the BVP is given by 
\[
\hat Y_k=Q_pJ_p^k(K^*K+E^*E)^{-1}K^*L,
\]
or, equivalently,
\[
\hat Y_k=\left[\begin{array}{ccc}1&1&0\\0&1&0\\0&0&0\\0&0&1\\0&0&1\end{array}\right]\left[\begin{array}{ccc}\frac{1}{2^k}&0&0\\0&0&0\\0&0&\frac{1}{4^k}\end{array}\right](K^*K+E^*E)^{-1}K^*L.
\]
Let
\[
E=\left[\begin{array}{ccc}0&\theta&0\\0&0&0\\0&0&0\\0&0&0\\0&0&0\end{array}\right]
\]
and $\theta=0.00001$. Then an optimal solution for the system is
%\[
%K^*K+E^*E=\left[\begin{array}{ccccc}1&1&0\\1&1+10^{-12}&0&0&0\\0&0&1&0&0\\0&0&0&0&1\end{array}\right]
%\]
%and 
\[
\hat Y_k=\left[\begin{array}{c}0\\0\\0\\\frac{11}{14\cdot 4^k}\\\frac{11}{14\cdot 4^k}\end{array}\right].
\]
\textit{Example 4.2.}\\\\
We assume the system \eqref{eq1} as in example 4.1 but with different boundary conditions. Let
\[
A_1=\left[\begin{array}{ccccc}1&0&0&0&0\\0&1&0&0&0\\0&0&1&0&0\\0&0&0&0&1\end{array}\right],
\]
\[
A_2=\left[\begin{array}{ccccc}1&0&0&0&1\end{array}\right],
\]
\[
B_1=\left[\begin{array}{c}0\\0\\0\\36\end{array}\right]
\]
and
\[
B_2=24.
\]
Or, equivalently,
\[
L=\left[\begin{array}{c}0\\0\\0\\36\\24\end{array}\right].
\]
With easy computations as in example 4.1, we get
\[
K=\left[\begin{array}{ccc} 36&36&0\\0&36&0\\0&0&36\\0&0&36\\358&35&24 \end{array}\right].
\]
We observe that $L \notin colspan K$, $r_1+r_2=5>3=p$, and $rank K=3=p$, i.e. the BVP is non-consistent. Thus from Theorem 4.1 and (a), an optimal solution of the BVP (3)-(4) is given by 
\[
\hat Y_k=Q_pJ_p^k(K^*K)^{-1}K^*L.
\]
or, equivalently,
\[
\hat Y_k=\left[\begin{array}{ccc}1&1&0\\0&1&0\\0&0&0\\0&0&1\\0&0&1\end{array}\right]\left[\begin{array}{ccc}\frac{1}{2^k}&0&0\\0&0&0\\0&0&\frac{1}{4^k}\end{array}\right]\left[\begin{array}{c} 0.0349\\-0.0166\\ 0.5006 \end{array}\right],
\]
or, equivalently,
\[
\hat Y_k=\left[\begin{array}{c}0.0349\cdot \frac{1}{2^k}\\0\\0\\0.5006\cdot\frac{1}{4^k}\\0.5006\cdot\frac{1}{4^k}\end{array}\right].
\]

\section*{Conclusions}
In this article we provided optimal solutions for a linear generalized discrete time system with boundary conditions whose coefficients are square constant matrices, its pencil is regular and the leading coefficient singular. 

%\subsection*{Acknowledgments} I would like to express my
%Our sincere gratitude to Professor G.I. Kalogeropoulos and Dr. I.K. Dassios for their helpful and
%fruitful discussions that improved this article.  

\end{document}